\documentclass[letterpaper, 10 pt, conference]{ieeeconf}  

\IEEEoverridecommandlockouts                              

\title{\LARGE \bf 
Decentralized Cooperative Merging of Platoons of Connected and Automated Vehicles at Highway On-Ramps
}
\author{Sharmila Devi Kumaravel, Andreas A. Malikopoulos,  Ramakalyan Ayyagari 
\thanks{This research was supported in part by ARPAE's NEXTCAR program under the award number DE- AR0000796 and by the Delaware Energy Institute (DEI). This support is gratefully acknowledged.}%
\thanks{Sharmila Devi Kumaravel and Ramakalyan Ayyagari are with Department of Instrumentation and Control Engineering, National Institute of Technology, Tiruchirappalli, India {\tt\small info2sj@gmail.com}; {\tt\small rkalyn@gmail.com}}%
\thanks{Andreas A. Malikopoulos is with the Department of Mechanical Engineering, University of Delaware, Newark, DE 19716 USA {\tt\small andreas@udel.edu}}%
}

\usepackage{multirow}
\usepackage{cite}  
\usepackage{amsmath,amssymb,amsfonts}
\usepackage{empheq}
\usepackage{ underscore}
\usepackage[linesnumbered,ruled,vlined]{algorithm2e}   \usepackage{multicol}
\usepackage{amssymb}
\usepackage{mathrsfs}
\usepackage{float,graphicx}
\usepackage{mathtools}
\usepackage{lscape}
\usepackage{setspace}
\def\mathclap#1{\text{\hbox to 0pt{\hss$\mathsurround=0pt#1$\hss}}}
\usepackage[mathlines]{lineno} 
\newcommand{\qedsymbol}{\rule{2mm}{2mm}}
\SetKwComment{Comment}{$\triangleright$\ }{}
\begin{document}
\maketitle
\thispagestyle{empty}
\pagestyle{empty}

\begin{abstract}

In this paper,  we present an optimization framework for cooperative merging of platoons of connected and automated vehicles at highway on-ramps. The framework includes (1) an optimal scheduling algorithm, through which, each platoon derives the sequence and time to enter the highway safely, and (2) the formulation of an optimal control problem the solution of which allows each platoon to derive its optimal control input (acceleration/deceleration) in terms of fuel consumption. We evaluate the efficacy of the proposed optimization framework through VISSIM-MATLAB simulation environment. The proposed framework significantly reduces the crossing time and fuel consumption of platoons at the highway on-ramps compared to the baseline scenario where the vehicles on the minor road yield to the vehicles on the highway. 
 
\end{abstract}

\section{INTRODUCTION}
Highway on-ramp is one of the main sources of traffic bottlenecks as the capacity of the road section is lesser than the neighboring roads. 
Several research efforts in the literature have used  ramp  metering  to  maximize the capacity of the highway by regulating the inflow of vehicles from on-ramps.  Athan \cite{c1} formulated a vehicle merging problem and facilitated the safe merging of vehicles without severe acceleration and deceleration except in emergency conditions. An overview of various research efforts in the literature that regulated  flow of vehicles  using ramp metering can be found in \cite{c2}. In these approaches, there is no communication between the merging vehicles and extreme care is expected from the drivers while  merging from on-ramps. The advances in connected vehicle technologies and automated driving can avoid human errors and offer better opportunities to  improve the mobility of commuters \cite{c3}.  

\par Cooperative merging control of connected and automated vehicles (CAVs)  can reduce congestion and provide significant benefits in terms of road capacity,  fuel consumption, travel time, and vehicle emissions \cite{c4}. Merging control algorithms for cooperative vehicles were presented in \cite{c5, c6} to create appropriate gap between vehicles from on-ramps to enter the motorway. In another effort, an optimal control framework for cooperative merging was presented to minimize fuel consumption \cite{c7}. A detailed summary  of various approaches for coordination of CAVs at highway on-ramps can be found in \cite{c4}. 

\par  Recently, an optimal control framework was developed to  minimize travel time and fuel consumption of vehicles \cite{c8}. The research effort in \cite{c9} formulated an optimal control problem to minimize the engine effort and passenger discomfort. Some research efforts facilitated safe merging  using optimization techniques \cite{c10}, reinforcement learning \cite{c11} and multi-agent Q-learning \cite{c12}.    The majority of research efforts in the literature have developed merging control algorithms without determining optimal merging sequence.  
 Several effective solutions to compute optimal sequence  can be achieved through scheduling theory \cite{c15}. Many scheduling theory-based control algorithms have been proposed in the literature to optimally coordinate the vehicles at various transportation scenarios \cite{c16,c17,c18,c19,c20,c21}.   
\par In our proposed framework, 
We model the  problem of cooperative merging at highway on-ramps as a job-shop scheduling problem, and present an optimal scheduling algorithm through which each platoon derives the sequence and time to enter the highway safely. Then, we present an optimal control problem where each platoon derives the optimal control input (acceleration/deceleration) in terms of fuel consumption of the platoons at the highway. 


\par The remainder of the paper is organized as follows. In Section II, we present the problem formulation and the optimal scheduling algorithm. In Section III, we provide an analytical solution for the optimal control problem. In Section IV, we evaluate the efficacy of the proposed framework through simulations, and present simulation results. We conclude with a discussion in Section V.

\section{PROBLEM FORMULATION}

\begin{figure*}[htbp!]
	\centering
	\includegraphics[width=4in]{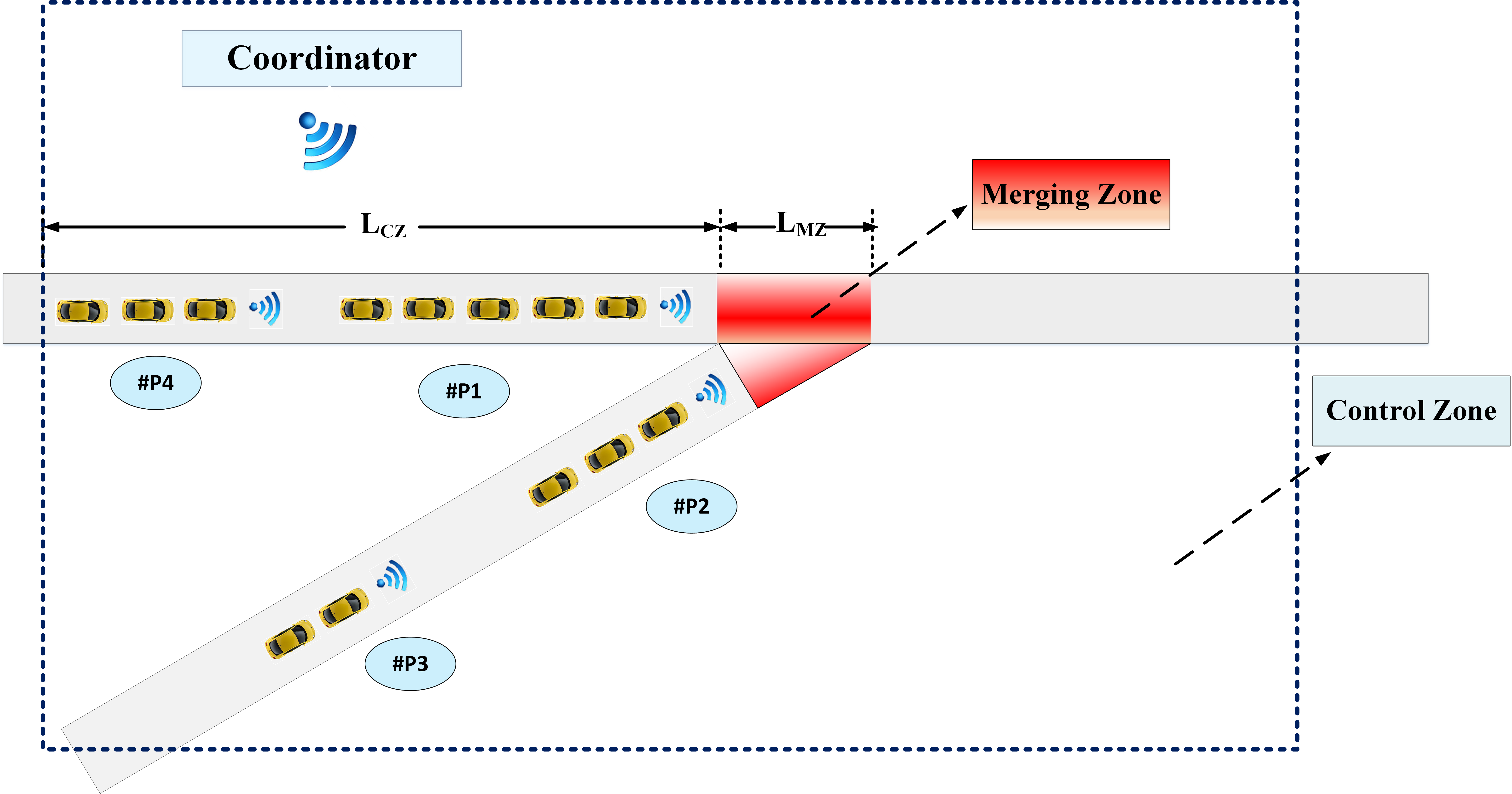}
	\caption{Platoons of connected and automated vehicles at highway on-ramps.}
	\label{platoonRamp}
\end{figure*}

We consider a traffic scenario where a secondary road with single-lane merges onto another single-lane main road  and P1, P2, P3 and P4 are the platoons of CAVs at merging roadways (Fig. \ref{platoonRamp}).
The region where the secondary road merges onto the main road is called the \emph{merging zone} (MZ) and has a length $L_{MZ}$. The merging zone is a conflict area where a potential lateral collision could occur. There is also a \emph{control zone} (CZ) inside of which the vehicles can communicate with each other. The length of the control zone is  $L_{CZ}$. A coordinator stores information about the geometry of the merging roadways and broadcasts this information back to the vehicles. Note that the coordinator does not take part in the decision-making process. In our modeling framework, we impose the following assumptions:\\\\
\textbf{Assumption 1:} 
The communication between the coordinator and CAVs occurs without any errors and delays.\\\\
\textbf{Assumption 2:} 
The platoons of CAVs enters the control zone as stable platoons, i.e., all the vehicles in the platoon move at a consensual speed and maintain the desired space between vehicles \cite{c22}. \\\\
\textbf{Assumption 3:} The value of $L_{CZ}$ should be sufficiently
large so that a platoon can reach the speed limit before entering the merging zone. \\\\
The first assumption may be strong, but it is relatively
straightforward to relax it as long as the measurement noise and delays are bounded. We impose the second assumption since our main focus is to facilitate the safe merging of CAVs entering as platoons at highway on-ramps rather than the stability and formation of the platoons. However, future research should relax this assumption and study the implications of the proposed solution on the stability and formation  of platoons. The third assumption is to enable the slower platoons with an initial speed less than the speed limit to reach the speed limit before they reach the merging zone. 

\subsection{Vehicle Model and Constraints}
Let $\displaystyle N(t) \in  \mathbb{N}$ be the number of platoons of CAVs arriving at the control zone at time $\displaystyle t \in \mathbb{R}^+$. Each platoon receives a unique identification number $\displaystyle k \in \mathbb{N}$ from the coordinator at the time of entry into control zone. The queue of platoons of CAVs in the control zone is denoted by $\mathcal N(t)=\{1,\ldots, N(t)\}$.  Let ${C}_k=\{1,\ldots,n_k\}$, $n_k\in \mathbb{N}$, be the number of vehicles in each platoon $\displaystyle k \in \mathcal N(t)$. We assume that each vehicle $\displaystyle j \in C_k$ is governed by the second order dynamics
\begin{align}
\dot{p}_j=v_j(t),\nonumber  \\ 
\dot{v}_j=u_j(t), 
\label{eq1}
\end{align}
where $\displaystyle p_j(t) \in \mathcal{P}_j$, $\displaystyle v_j(t) \in \mathcal{V}_j$,  $\displaystyle u_j(t) \in \mathcal{U}_j$ denote position, speed, and,  acceleration/deceleration. 
\par Let $\displaystyle t_{j}^{0}$ and $t_{j}^{m}$ be the time at which the leader vehicle $\displaystyle j \in C_k$ enters the control zone and merging zone, respectively. Let $\displaystyle x_j(t)=[p_j(t) \quad v_j(t)]^T$ denote the state of each vehicle $\displaystyle j \in C_k$ with ${x_j}^0=[p_j^0 \quad v_j^0]^T$ as initial state, where $\displaystyle  p_j^0=p_j(t_j^0)=0$, taking values in the state space $\displaystyle \mathcal{X}_j=\mathcal{P}_j \times \mathcal{V}_j$.   The control input and speed of each vehicle $\displaystyle j \in C_k$ is bounded with  the following constraints
\begin{gather}
u_{\min} \leq u_j(t)  \leq u_{\max},  \label{eq2} \\
0 \leq v_{\min} \leq v_j(t) \leq v_{\max}, \label{eq3} 
\end{gather}
where $\displaystyle u_{\min}, u_{\max}$ are the minimum and maximum control inputs  and $\displaystyle v_{\min}, v_{\max}$  are the minimum and maximum speed limits, respectively.
\subsection{Optimal Scheduling Algorithm for Platoons Crossing the Merging Zone}
 In this section, we apply scheduling theory to derive a control algorithm that provides the optimal schedule for each platoon to enter the merging zone. We consider a job-shop scheduling problem to find the  optimal sequence and  schedule  for the platoons at merging roadways. 
 A job-shop scheduling problem addresses the scheduling of multiple jobs on a single machine or several machines. It  optimally allocates the jobs on the machines with an objective to improve the performance criteria. We model the problem of cooperative merging of platoons at highway on-ramps as a job-shop scheduling problem. We consider merging roadways as a single machine and the platoons of CAVs as jobs.

  \par Let $t_k^{\text{in}}$ be the time that a platoon $\displaystyle k \in \mathcal{N}(t)$ takes to reach the merging zone. We consider the speed at which the platoon enters the control zone as initial speed of the platoon. We consider two cases to compute the time at which each platoon reaches the merging zone. In Case 1, we consider that the initial speed of platoons is the same as the speed limit.  In Case 2, we consider that the initial speed of platoons is less than the speed limit.  Let $\displaystyle v_{k}^0= v_k(t_k^0)$ be the initial speed of the platoon $\displaystyle k \in C_k$. 
\newpage
\noindent \textbf{Case 1: $\displaystyle v_{k}^0= v_{\max}$}\\
We compute $\displaystyle {t_{k}^{\text{in}}}$   using  $\displaystyle v_{\max}$, i.e., 
\begin{align}
t_{k}^{in}=\frac{L_{CZ}}{v_{\max}}.\label{eq6}
\end{align}
\textbf{Case 2: $\displaystyle v_{k}^0< v_{\max}$} \\
Let $\displaystyle t_{k}^a$ be the time taken by the platoon $k$ to accelerate to the speed limit. Then, 
\begin{align}
t_{k}^a=\frac{v_{\max}-v_{k}^0}{u_{\max}}.\label{eq7}
\end{align}
Let  $\displaystyle d_{k}^a$ be the distance traveled during acceleration, then,
\begin{align}
d_{k}^a=\frac{{({v_{\max}})^2}-(v_{k}^0)^2}{2u_{\max}}, \label{eq8}
\end{align}
and thus,
\begin{align}
t_{k}^{in}= t_{k}^a+ \frac{L_{CZ}-d_{k}^a}{v_{\max}}.\label{eq9}
\end{align}
Let $t_k^{\text{out}}$ be the time that a platoon $k \in \mathcal N(t)$ takes to exit the merging zone. Then, 
 \begin{align}
   t_{k}^{\text{out}}=\frac{L_{MZ}}{v_{\max}} + (|C_k|-1) \times t_k^{h}+t_{g}, \label{eq10}
  \end{align} 
where $|C_k|$ is the number of CAVs in platoon $\displaystyle k \in \mathcal N(t)$,  $\displaystyle t_{g}$ is the safe time gap provided to ensure safety of platoons entering and leaving the merging zone, and $\displaystyle t_{k}^h$ is the headway between consecutive CAVs in the platoon $k$.\\\\
\textbf{Definition 1:} The completion time $t_k^c$ is the time taken by the platoon to completely exit the merging zone.The completion time $t_k^c$ of platoon $\displaystyle k \in \mathcal{N}(t)$  is 
 \begin{align} 
t_k^c \coloneqq t_k^{\text{in}}+ t_k^{\text{out}}.\label{eq11}
 \end{align}
\textbf{Definition 2:} The total completion time $T_C$ of platoons to cross the merging zone is  
 \begin{align}
 T_C\coloneqq \sum\limits_{k=1}^N t_{k}^c.\label{eq4}
 \end{align}
\textbf{Definition 3:} Let $w_k$ be the weight assigned to each platoon  $k \in \mathcal N(t)$  based on the priority.  The total weighted completion time $T_{WC}$  of platoons is  
 \begin{align}
 T_{WC}\coloneqq \sum\limits_{k=1}^N w_kt_{k}^c.\label{eq5}
 \end{align}
 We model the problem of cooperative merging of platoons of CAVs at highway on-ramps as  $1||T_{WC}$ job-shop scheduling problem \cite{c23}, where 1 represents a single machine, and $T_{WC}$ represents the total weighted completion time. A schedule is optimal if it minimizes $T_{WC}$. The objective of the modeled $1||T_{WC}$ scheduling problem is to minimize the crossing time, i.e., the total weighted completion time of platoons at highway on-ramps. Each platoon leader solves $1||T_{WC}$ problem to derive the time to enter the merging zone with the objective to minimize the crossing time of platoons at the merging zone, i.e., total weighted completion time of platoons. The platoons entering the merging zone from the main roadway are given higher priority than platoons entering from the secondary road. The higher priority to platoons in main roadway allow them  to enter the merging zone before the platoons entering from the secondary road.  The weights $w_k$ are assigned based on the priority of the platoon. We aim to avoid the deceleration of vehicles in the main roadway. Hence, $w_k$ for platoons entering the main roadway is assigned greater value than the the platoons entering from the secondary road. 
 Each platoon leader communicates with the coordinator to acquire the information about the geometry and topology of the merging roadways. In addition, the platoon leaders broadcast their attributes including the link number, lane number, unique identification number, the speed at which they enter the control zone, maximum acceleration, maximum speed limit, number of followers in the platoon, safe time gap, headway and speed limit to other platoon leaders and coordinator.
 \vspace{-1mm}
  \par The optimal scheduling algorithm running in each platoon leader uses the information of initial speed, maximum acceleration, and maximum speed limit  to compute the time that a platoon arrives at the merging zone. Then, it uses the information of the length of the merging zone, number of followers in the platoon, safe time gap, headway of each platoon, and speed limit to compute the time that the platoon exits the merging zone. Next, it computes the completion time  of each platoon based on the weights assigned to the platoons. Then, it calculates the total weighted completion time of platoons. Further, the algorithm  arranges the platoons in the non-decreasing order of weighted completion time. Finally, the algorithm computes the time at which each platoon has to enter the merging zone. The proposed scheduling algorithm yields an optimal
schedule that minimizes the crossing time of platoons.\\
  \begin{algorithm*}
 	\SetAlgoLined
 \caption{ Computation of optimal sequence and schedule of platoons}
\begin{multicols}{2}
	\KwData{$v_{k}^0$, $w_k$, $|C_k|$, $t_{k}^h$ of each platoon $k$, current time $t$, $t_{g}$, $L_{CZ}$, $L_{MZ}$, $v_{\max}$, $u_{\max}$}
	\KwResult{$optimalSequence$ and schedule $t_i^m$ for platoons}
	Initialize variable $t_{last}^l \gets 0$\\
    \Comment{computation of total weighted completion time of platoons}
 	\For {$k=1~ to~ N$}{
 	\uIf {$v_{\max}=v_{k}^0$}{
 	 $t_{k}^{\text{in}} \gets L_{CZ}/v_{\max}$
 }
 	\ElseIf {$v_{\max}<v_{k}^0$}{
 	 $t_{k}^{a} \gets (v_{\max}-v_{k}^0)/ u_{\max}$\\
 	$d_{k}^a \gets (v_{\max})^2-(v_{k}^0)^2/ 2u_{\max}$\\
 	$t_k^{in} \gets t_{k}^a+(L_{CZ}-d_{k}^a)/v_{\max}$
 }	
	    
 }
 	\For {$k=1~to~N$}{
 	 $t_{k}^{\text{out}} \gets L_{MZ}/v_{\max}+(|C_k|-1)*t_k^{h}+t_{g}$
}
 	\For {$k=1~to~N$}{
 	 $t_{k}^c \gets t_{k}^{in} +t_{k}^{\text{out}}$\\
 	 $w_k^c \gets t_{k}^c/w_k$
 }
 	 Arrange platoons in non-decreasing order of $w_k^c$ in an array $optimalSequence$\\	
 	 \Comment{Computation of time of entry of CAVs in Platoons}
     $i=optimalSequence[1]$\\
 	 $t_{i}^{m} \gets t + t_{i}^{\text{in}}$\\
 	 $t_{i}^l\gets t+ t_{i}^{c}$ \\	
 	\If{$t_i^m < t_{last}^l$}{
 	 $t_i^m \gets t_{last}^l$\\	
 	 $t_{i}^l \gets t_{last}^l+t_{i}^{\text{out}}$
 	 
 }
 	\For {$j=2~ to~ N$}{
 	 $r=optimalSequence[j]$\\
 	$t_{r}^{m} \gets t_{r-1}^{l}$\\
 	 $t_{r}^l \gets t_{r}^{m}$ + $t_{r}^{\text{out}}$\\
 	\If {$j==N$}{
 	 $last=optimalSequence[j]$\\
 	 $t_{last}^l \gets t_{j}^l$\\
 }
}
\end{multicols}
\end{algorithm*}
\noindent \textbf{Theorem 1:} \emph{A schedule $\displaystyle \Lambda^*$  is optimal in minimizing the crossing time of the platoons, i.e., the total weighted completion time if and only if it places the platoons in the order of non-decreasing weighted completion time $\displaystyle w_{k}^c$. } \\\\
\textbf{Proof:}
Let $p$, $q$, and $r \in \mathcal{N}(t)$ be the platoons in the control zone arriving at time $\displaystyle t$ and $\frac{t_r^c}{w_r} < \frac{t_p^c}{w_p} < \frac{t_q^c}{w_q} $. Let $\Lambda$ be the schedule in which platoons are not arranged in the order of non-decreasing  weighted completion times $\displaystyle w_{k}^c$ of platoons. Then in schedule $\Lambda$,  platoon $q$ precedes platoon $p$ and $\frac{t_q^c}{w_q} > \frac{t_p^c}{w_p}$. Let  $\Lambda^*$ be the another schedule in which platoon $p$ precedes platoon $q$.  We know that $t_p^c$, $t_q^c$, and   $t_r^c$ are the completion time of platoons $p$, $q$, and $r$, respectively. Let $w_p$,  $w_q$, and $w_r$ be the weights assigned to platoons $p$, $q$, and $r$, respectively.\\
Let $T_{WC}(\Lambda)$ and $T_{WC}(\Lambda^*)$ be the total weighted completion time of the platoons in the schedule $\Lambda$ and $\Lambda^*$, respectively. The total weighted completion time of schedules $(\Lambda)$ and $(\Lambda^*)$ are
\begin{eqnarray}
T_{WC}(\Lambda)&=& w_r(t+t_r^c) + w_q(t+t_r^c+t_q^c) \nonumber\\ && + w_p(t+t_r^c+t_q^c+t_p^c),\label{eq12}
\end{eqnarray}
\begin{eqnarray}
T_{WC}(\Lambda^*)&=&w_r(t+t_r^c) + w_p(t+t_r^c+t_p^c) \nonumber\\ && + w_q(t+t_r^c+t_p^c+t_q^c).\label{eq13}
\end{eqnarray}
\begin{equation}
   T_{WC}(\Lambda^*)- T_{WC}(\Lambda)=w_qt_p^c-w_pt_q^c < 0\label{eq14}. 
\end{equation}
This contradicts the optimality of schedule $\Lambda$.  Hence, schedule $\Lambda^*$ in which platoons  arranged in the order of non-decreasing weighted completion time $\displaystyle w_{k}^c$  is optimal.
\hspace*{80mm}\qedsymbol\\
 \subsection{Optimal Control Problem}
Each platoon leader computes the time to enter the merging zone using the optimal scheduling algorithm. Then, each leader derives the optimal control input to enter the merging zone in a schedule specified by the scheduling algorithm.  The position of platoons in the queue designates the time that a platoon enters the merging zone. Here, we have two cases: 1) If  $t_k^m=t_k^{in}$, then the leader solves a time optimal control problem, and 2) if $t_k^m>t_k^{in}$, then the leader solves an energy optimal control problem to derive the optimal control input. The leader broadcasts the derived schedule and optimal control input to the following vehicles until the last follower has completely crossed the merging zone.
\subsubsection{Time Optimal Control Problem}
For each leader $\displaystyle  j \in C_k$, we formulate the following time optimal control problem,  
\begin{align}
\min_{u_{j}\in \mathcal{U}_j}J_1^j(u_j(t))=\int_{t_j^0}^{t_j^m}dt= t_{j}^{m}-t_{j}^{0},\label{eq15}
\end{align}
subject to: (\ref{eq1}), (\ref{eq2}), (\ref{eq3}), $\displaystyle p_j(t_j^0)=0$, $\displaystyle p_j(t_{j}^{m})=L_{CZ}$, \\ and given  $\displaystyle t_j^0, v_{j}^0, t_{j}^{m}$,\\ where $t_j^0$ and $t_{j}^{m}$ is the time at which the leader enters the control zone and merging zone, respectively.
\subsubsection{Energy Optimal Control Problem}
For each leader $\displaystyle j \in C_k$, we formulate the following energy optimal control problem, 
\begin{align}
\min_{u_{j}\in \mathcal{U}_j} J_2^j(u_j(t))=\frac{1}{2}\int_{t_j^0}^{t_j^m} u_j^2(t) \; dt, \label{eq16} 
\end{align}
subject to: (\ref{eq1}), (\ref{eq2}), (\ref{eq3}),  $\displaystyle p_j(t_j^0)=0$,  $\displaystyle p_j(t_{j}^{m})=L_{CZ}$ and given  $\displaystyle t_j^0, v_{j}^0, t_{j}^{m}$,
\\\\
where $t_j^0$ and $t_{j}^{m}$ is the time at which the leader enters the control zone and merging zone, respectively.
\vspace{-1mm}
\section{Analytical Solution}
 In this section, we present the closed-form analytical solutions for the time (\ref{eq15}) and energy (\ref{eq16}) optimal control problems for each leader $\displaystyle j \in C_k$.

\subsection{Analytical Solution of the Time Optimal control problem}
We derive the Hamiltonian function for each leader $j\in C_k$ with the state and control constraints as follows
\begin{gather}
H_j(t,p_j(t),v_j(t),u_j(t))=1 + {\lambda_j^p}{ v_j(t)}+ \lambda_j^v u_j(t) \nonumber \\
+ \mu_j^a (u_j(t)-u_{\max}) + \mu_j^b (u_{\min}-u_j(t))\nonumber \\+\mu_j^c (v_j(t)-v_{\max}) + \mu_j^d (v_{\min}-v_j(t)), \label{eq17} 
\end{gather} 
where $\displaystyle {\lambda_j^p}$ and $\displaystyle {\lambda_j^v}$ are costates and  $\displaystyle {\mu_j^a}$,  $\displaystyle {\mu_j^b}$, $\displaystyle {\mu_j^c}$, and $\displaystyle {\mu_j^d}$ are Lagrange multipliers. We derive optimal control input for platoons with initial speed either equal to speed limit or less than speed limit while entering the control zone.\\
\textbf{Case 1:}  $v_j^0\le v_{\max}$\\
In this case, the optimal control input is
\begin{equation}
u_j^*(t)=\begin{cases}
u_{\max}, & \text{if $t_{j}^{0} \leq t \leq t_{j}^{0}+ t_{j}^a$},\\
0, & \text{if $t_{j}^{0}+ t_{j}^a \leq t \leq t_{j}^{m}$}.
\end{cases}\label{eq18} 
\end{equation}
Substituting (\ref{eq18}) in (\ref{eq1}), we derive the optimal position and velocity
 \begin{empheq}[left={\empheqlbrace}]{align}
&p_j^*(t)=\frac{1}{2}u_jt^2+b_jt+c_j,&\text{if}~ t_{j}^{0} \leq t \leq t_{j}^{0}+ t_{j}^a \label{eq19}\\
&v_j^*(t)=u_jt+b_j,&\text{if}~ t_{j}^{0} \leq t \leq t_{j}^{0}+ t_{j}^a \label{eq20}\\
&p_j^*(t)=v_{\max}t+d_j,&\text{if}~t_{j}^{0}+ t_{j}^a \leq t \leq t_{j}^{m}\label{eq21}\\
&v_j^*(t)=v_{\max},&\text{if}~t_{j}^{0}+ t_{j}^a \leq t \leq t_{j}^{m}\label{eq22}
\end{empheq}
where $ b_j$, $c_j$, and $d_j$ are integration constants computed  using the initial and final conditions in (\ref{eq15}).\\
\textbf{Case 2:}  $v_j^0= v_{\max}$\\
The optimal control input based on \cite{c24} is 
\begin{equation}
u_j^*(t)= 0,\quad{ t \in [ t_{j}^{0},\quad t_{j}^{m}]}.
\label{eq23} 
\end{equation}
Substituting (\ref{eq23}) in (\ref{eq1}), we derive the optimal position and velocity,
\begin{align}
p_j^*(t) &= v_{\max}t+d_j,\label{eq24}\quad t \in [t_{j}^{0},~ t_{j}^{m}],\\
v_j^*(t) &= v_{\max},\quad \label{eq25} t \in [t_{j}^{0},\quad t_{j}^{m}],
\end{align}
where $ d_j$ is integration constant computed  using the initial and final conditions in (\ref{eq15}). 
\subsection{Analytical Solution of the Energy Optimal Control Problem}
  For this problem, the Hamiltonian function for each leader $j\in C_k$ is
\begin{gather}
H_j(t,p_j(t),v_j(t),u_j(t))=\frac{1}{2} u_j^2(t) + {\lambda_j^p}{ v_j(t)}+ \lambda_j^v u_j(t) \nonumber \\
+ \mu_j^a (u_j(t)-u_{j,\max}) + \mu_j^b (u_{j, \min}-u_j(t))\nonumber \\+\mu_j^c (v_j(t)-v_{j, \max}) 
+ \mu_j^d (v_{j, \min}-v_j(t)),\label{eq26}
\end{gather}
where $\displaystyle {\lambda_j^p}$ and $\displaystyle {\lambda_j^v}$ are costates, and  $\displaystyle {\mu_j^a}$,  $\displaystyle {\mu_j^b}$, $\displaystyle {\mu_j^c}$, and $\displaystyle {\mu_j^d}$ are the Lagrange multipliers and $\mu_j^a = \mu_j^b = \mu_j^c = \mu_j^d = 0$ as the state and control constraints are not active. The optimal control input based on \cite{c25} is
\begin{gather}
u_j^*=a_jt+b_j,\quad\label{eq27}t \in [t_{j}^{0},\quad t_{j}^{m}].
\end{gather}
substituting (\ref{eq27}) in (\ref{eq1}), we get the optimal position and velocity,
\begin{align}
p_j^*(t)&=\frac{1}{6}a_jt^3+\frac{1}{2}b_jt^2+c_jt+d_j,\quad t \in [t_{j}^{0},\quad t_{j}^{m}], \label{eq28}\\
v_j^*(t)&=\frac{1}{2}a_jt^2+b_jt+c_j,\quad t \in [t_{j}^{0},\quad t_{j}^{m}],\label{eq29}
\end{align}
where $a_j,~ b_j,~ c_j$, and $d_j$ are integration constants computed  using initial and final conditions in (\ref{eq16}).
\section{Simulation Results and Discussions}
We model the traffic scenario where secondary road with single lane merges onto single lane  main road in  VISSIM 11.00 simulation environment. 
The size of platoons entering the main roadway varies from $1$ to $5$ vehicles. The size of platoons entering from the secondary road ranges from $1$ to $3$ vehicles. The length of the merging zone  is designed to be  $30~m$.  Based on Assumption 3, the length of the control zone is designed to be $150~m$ so that the platoons can reach the speed limit before entering the merging zone. The maximum speed limit is set as $25~m/s$.  The maximum acceleration  and minimum deceleration limit are $3~m/s^2$ and $-3~ m/s^2$, respectively. 
\par The optimal scheduling algorithm is implemented in MATLAB. 
The information of the platoons 
are collected in real-time through COM interface from the VISSIM simulation environment. Based on the collected information, each platoon leader calculates the time to enter the merging zone. Further, the optimal control input for each platoon to enter the merging zone is derived. The speed of platoons corresponding to the derived optimal control input is updated in real-time through COM interface in VISSIM. 
\par
\begin{figure}[htbp!]
	\centering
	\includegraphics[scale=0.5]{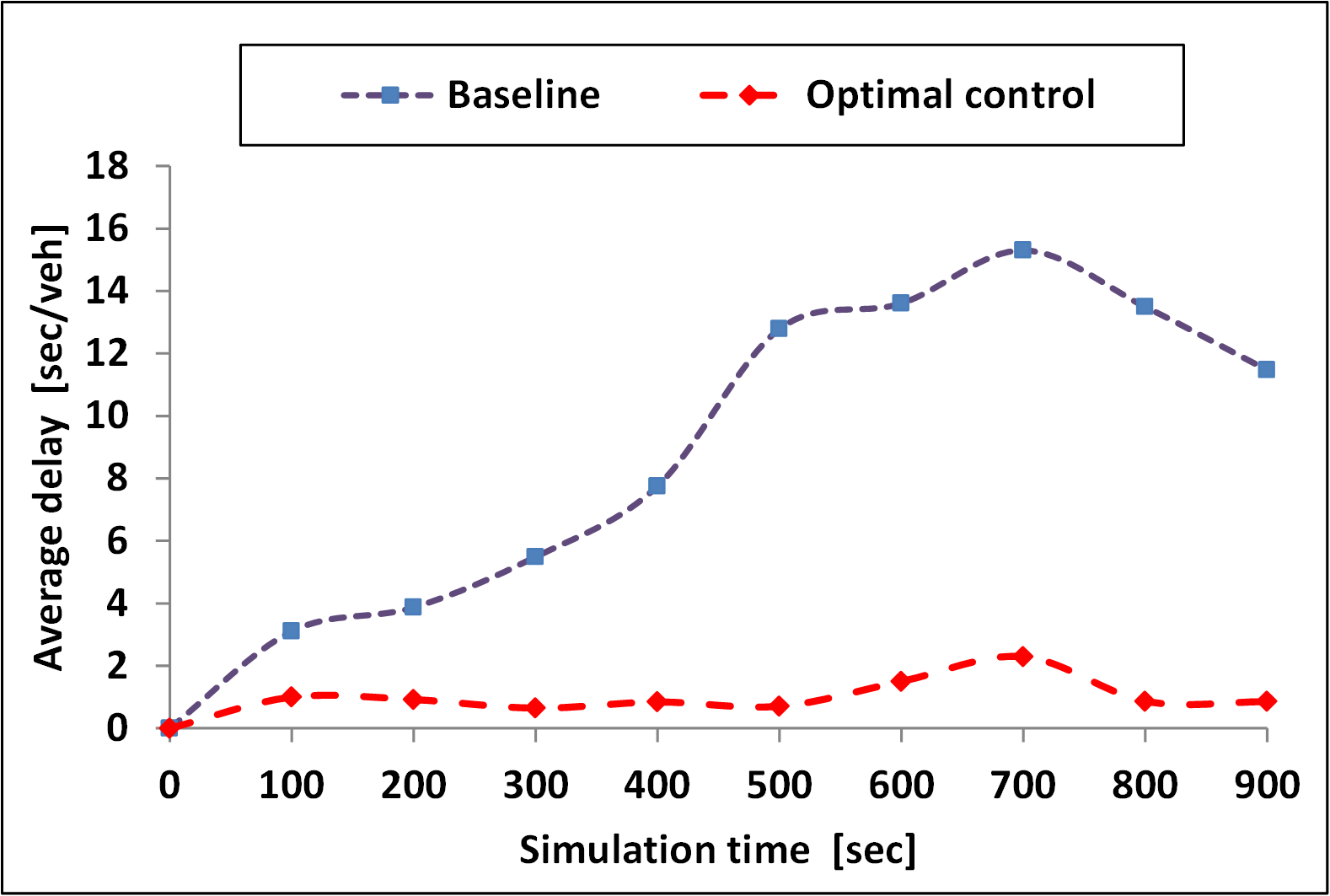}
	\caption{Average delay.}
	\label{delay}
\end{figure}

\begin{figure}[htbp!]
	\centering
	\includegraphics[scale=0.5]{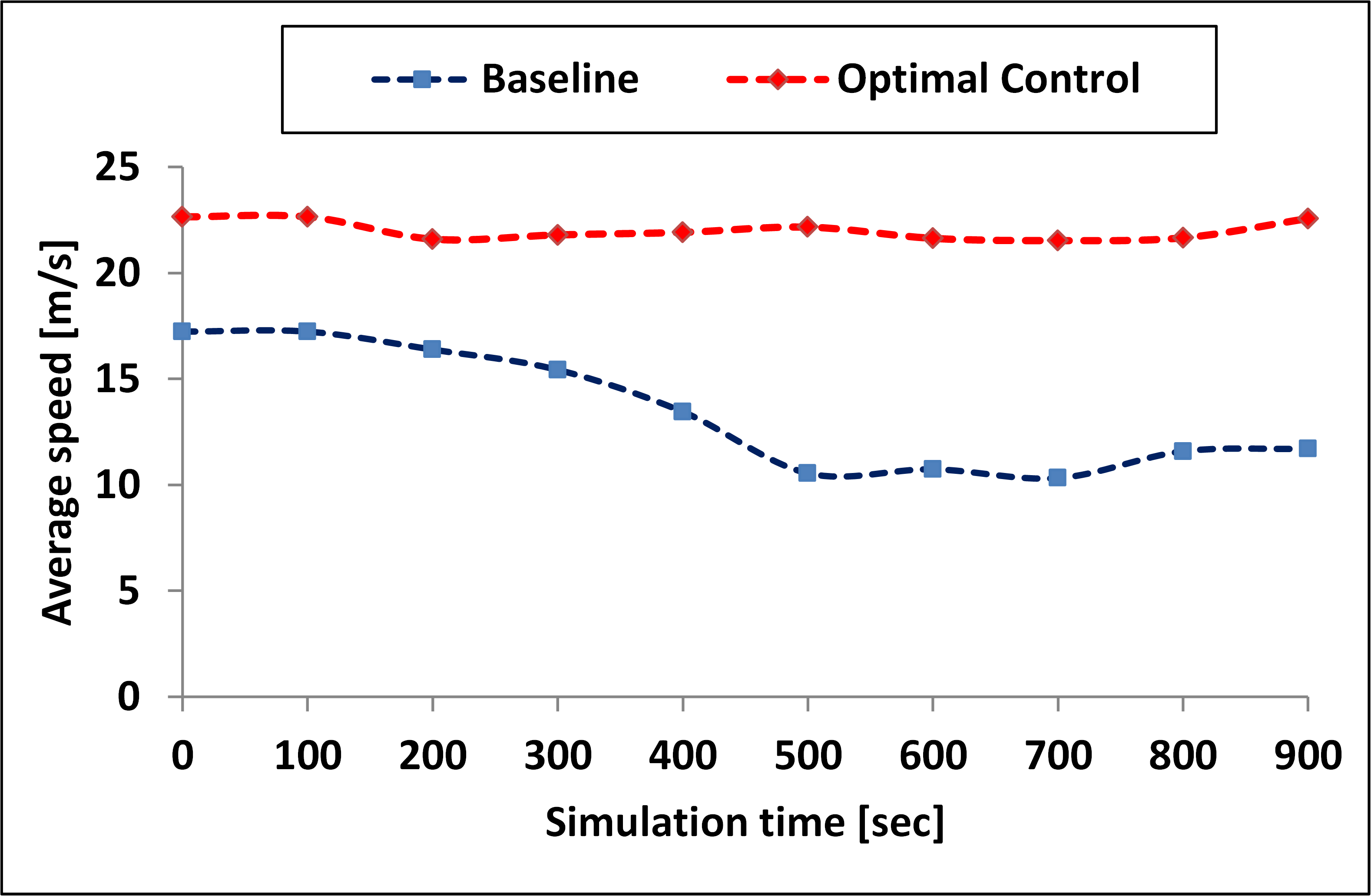}
	\caption{Average speed.}
	\label{avgspeed1}
\end{figure}

\begin{figure}[htbp!]
	\centering
	\includegraphics[scale=0.5]{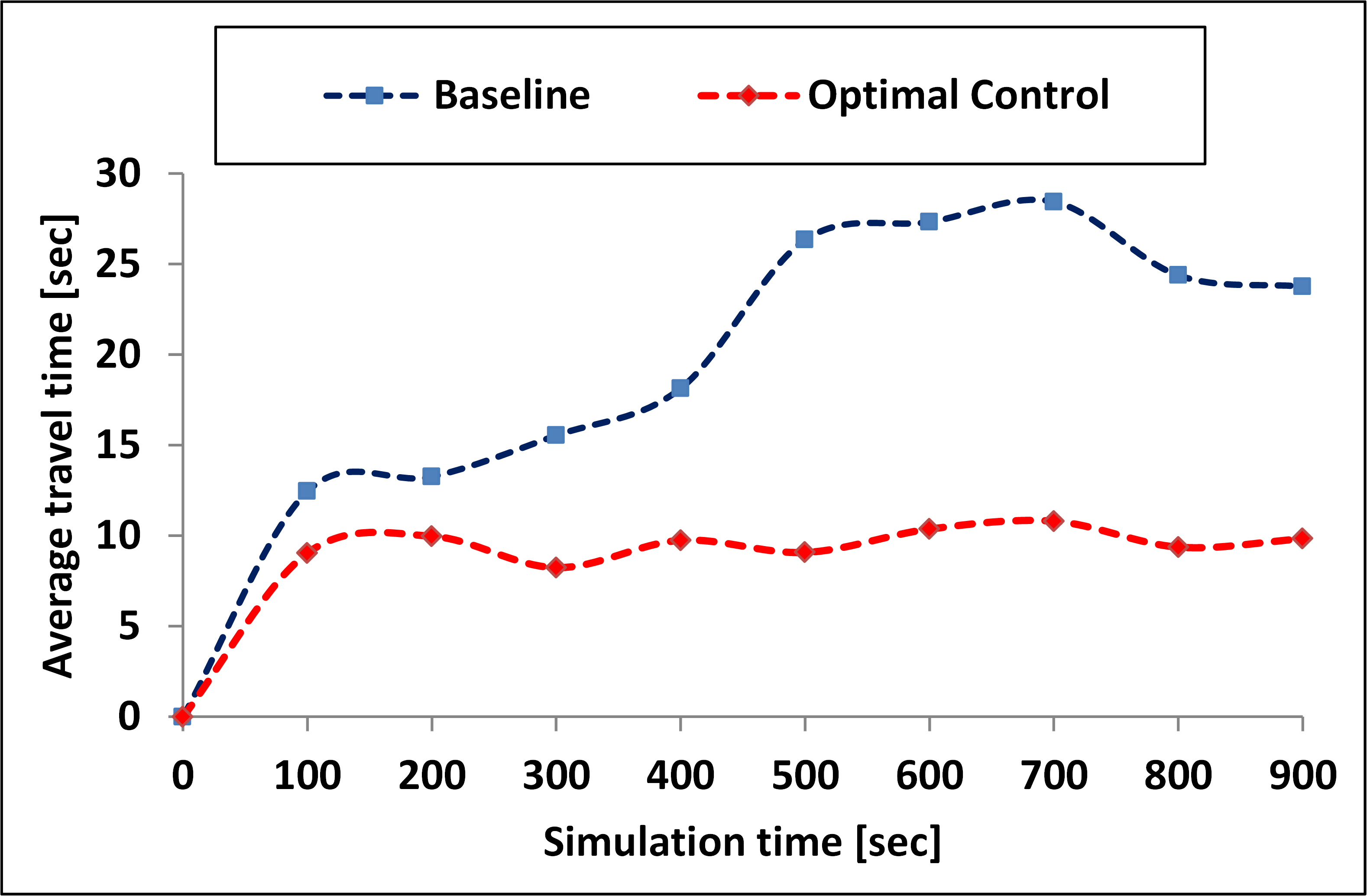}
	\caption{Average travel time.}
	\label{avgtraveltime}	
\end{figure}

\begin{figure}[htbp!]
	\centering
	\includegraphics[scale=0.5]{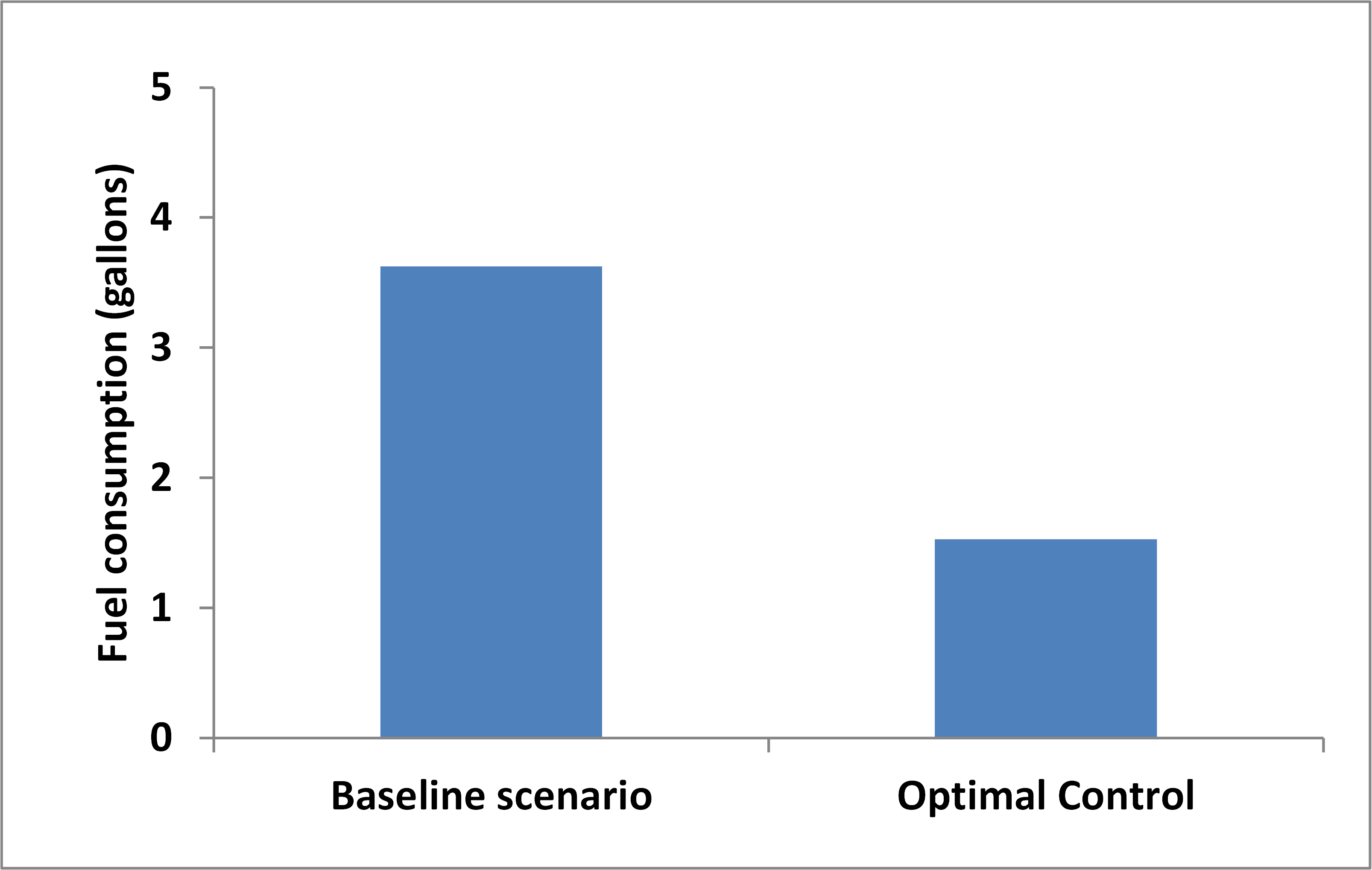}
	\caption{Fuel consumption.}
	\label{fuel}	
\end{figure}

 We compare the performance of the proposed framework with a baseline scenario where the vehicles from secondary road stop and yield for vehicles entering from the main road. We design the volume of vehicles entering from the main road as $1060~vph$. The volume of vehicles entering the secondary road is considered as $720~vph$. The simulation is run for $900$ seconds. The performance measures of the vehicles under baseline scenario and proposed optimal framework are  shown in Figs. \ref{delay},    \ref{avgspeed1}, \ref{avgtraveltime}, and \ref{fuel}. The proposed framework eliminated the stop-and-go driving behavior and resulted in negligible delay for platoons.  The proposed framework reduced average travel time and fuel consumption by 54.3\% and 57.8\%, respectively. The average delay is reduced by 88.92\% and the average speed of platoons is increased by 63.53\% compared to the baseline scenario. 

 \section{Concluding Remarks and Discussions}
In this paper, we presented an  optimization framework for cooperative merging of platoons of CAVs at highway on-ramps. We modeled the problem as a job-shop scheduling problem and presented an algorithm to derive the schedule of each platoon to enter the merging zone. Then, we provided the analytical solution for the optimal control problem to minimize fuel consumption and crossing time of platoons at highway on-ramps. We validated the proposed framework through simulations, and we noticed significant improvements in fuel consumption, while minimizing stop-and-go driving behavior.
 
Ongoing efforts consider lane changes for coordination of CAVs at multi-lane roads \cite{c26} and multiple scenarios \cite{c27}. In our framework, we assumed that vehicles form stable platoons before arriving  the control zone. Future research should focus on the stability and formation of the platoons along with cooperative merging  in mixed traffic environment  with  human-driven vehicles and CAVs.


\begin{thebibliography}{99}
		\bibitem{c1}
		M. Athans, ``A unified approach to the vehicle-merging problem,'' \emph{Transp. Res.},vol. 3, no. 1, pp. 123\textendash133, Apr. 1969.
	\bibitem{c2}
	M. Papageorgiou and A. Kotsialos, ``Freeway ramp metering: An overview,'' \emph{IEEE transactions on intelligent transportation systems}, vol. 3, no. 4, pp. 271\textendash281, 2002.
	\bibitem{c3}
	L. Zhao and A. A. Malikopoulos, ``Enhanced mobility with connectivity and automation: A review of shared autonomous vehicle systems,'' \emph{IEEE Intelligent Transportation Systems Magazine}, 2020.
		\bibitem{c4}
	J. Rios-Torres and A. A. Malikopoulos, ``A survey on the coordination of connected and automated vehicles at intersections and merging at highway on-ramps,'' \emph{IEEE Transactions on Intelligent Transportation Systems}, vol. 18, no. 5, pp. 1066\textendash1077, 2017.

	\bibitem{c5}
    R. Scarinci, B. Heydecker, and  A. Hegyi, ``Analysis of traffic performance of a merging assistant strategy using cooperative vehicles,'' \emph{IEEE Transactions on Intelligent Transportation Systems}, vol. 16, no. 4, pp. 2094\textendash2103, 2015.  
    
    
    \bibitem{c6}
   A. Mosebach, S. Röchner, and J. Lunze, ``Merging control of cooperative vehicles,'' \emph{IFAC-PapersOnLine}, vol. 49, no. 11, pp. 168\textendash174, 2016.
   \bibitem{c7}
   J. Rios-Torres and A. A. Malikopoulos, ``Automated and cooperative vehicle merging at highway on-ramps,'' \emph{IEEE Transactions on Intelligent Transportation Systems}, vol. 18, no. 4, pp. 780\textendash789, 2017.	 

	
	\bibitem{c8}
	W. Xiao  and C. G. Cassandras, ``Decentralized optimal merging control for connected and automated vehicles,'' \emph{ In 2019 American Control Conference (ACC)}, pp. 3315\textendash3320, 2019.
	\bibitem{c9}
      I. A. Ntousakis, I. K. Nikolos,  and M. Papageorgiou, ``Optimal vehicle trajectory planning in the context of cooperative merging on highways,'' \emph{Transportation research part C: emerging technologies}, vol. 71, pp.  464\textendash488, 2016.
      	\bibitem{c10}
		 G. Raravi, V. Shingde, K. Ramamritham, and J. Bharadia, ``Merge algorithms for intelligent vehicles,'' \emph{Proc. Next Gener. Des. Verif. Methodol. Distrib. Embed. Control Syst.}, pp. 51\textendash65, 2007.
	\bibitem{c11}
	P. Wang and C. Y. Chan, ``Autonomous ramp merge maneuver based on reinforcement learning with continuous action space,'' \emph{arXiv preprint arXiv:1803.09203}, 2018.
	\bibitem{c12}
	L. Schester and L. E. Ortiz, ``Longitudinal position control for highway on-ramp merging: A multi-agent approach to automated Driving,'' \emph{ In 2019 IEEE Intelligent Transportation Systems Conference} pp. 3461-3468, 2019.   
		
		\bibitem{c15}
		P. Brucker and P. Brucker, ``Scheduling algorithms,'' Vol. 3. Berlin: Springer, 2007.
	\bibitem{c16}
J. Wu, A. Abbas-Turki, and A. El Moudni, ``Cooperative driving: an ant colony system for autonomous intersection management,'' \emph{Applied Intelligence}, vol. 37, no. 2, pp. 207\textendash222, 2012.
		\bibitem{c17}
L. Bruni, A. Colombo, and D. Del Vecchio, ``Robust multi-agent collision avoidance through scheduling,'' \emph {In 52nd IEEE Conference on Decision and Control}, pp. 3944\textendash3950, 2013.
	\bibitem{c18}
	N. Murgovski, G. R. de Campos, and J. Sjöberg,   ``Convex modeling of conflict resolution at traffic intersections,'' \emph{In 2015 54th IEEE conference on decision and control}, pp. 4708\textendash4713, 2015.
	\bibitem{c19}
	G.R. De Campos, F. Della Rossa, and A. Colombo,  ``Optimal and least restrictive supervisory control: Safety verification methods for human-driven vehicles at traffic intersections, \emph{In 2015 54th IEEE Conference on Decision and Control}, pp. 1707\textendash1712, 2015.
	\bibitem{c20}
		B. Chalaki and A. A. Malikopoulos, ``Time-optimal coordination
for connected and automated vehicles at adjacent intersections,'' \emph{arXiv preprint arXiv:1911.04082}, 2020.
	\bibitem{c21}
 M. A. Guney and I. A. Raptis, ``Scheduling-Based Optimization for Motion Coordination of Autonomous Vehicles at Multilane Intersections,`` \emph{Journal of Robotics}, 2020.
  \bibitem{c22}
 W. S. Levine,  and M. Athans, ``On the optimal error regulation of a string of moving vehicles,'' \emph{IEEE Transactions on Automatic Control}, Vol. 11, no. 3, pp.355-361.  1966. 
  
	\bibitem{c23}
	W. E. Smith, ``Various optimizers for single‐stage production,'' \emph{Naval Research Logistics Quarterly}, Vol. 3, no. 1‐2, pp. 59\textendash66, 1956.
		\bibitem{c24}
		S. D. Kumaravel, A. A. Malikopoulos, and R. Ayyagari, ``Optimal coordination of platoons of connected and automated vehicles at signal-free intersections,'' \emph{arXiv preprint arXiv:2001.04866}, 2020.
			\bibitem{c25}
			A. A. Malikopoulos, C. G. Cassandras, and Y. Zhang, ``A decentralized energy-optimal control framework for connected automated vehicles at signal-free intersections,''\emph{Automatica}, vol. 93, pp. 244\textendash256, 2018.
			
			\bibitem{c26}
			A. A. Malikopoulos,  Beaver, L.E., and Chremos, I.V., ``Optimal Time Trajectory and Coordination for Connected and Automated Vehicles,''\emph{Automatica}, vol. 125, pp. 109469 2021.
			
			\bibitem{c27}
		B. Chalaki and A. A. Malikopoulos, ``Optimal Control of Connected and Automated Vehicles at Multiple Adjacent Intersections,'' \emph{IEEE Trans. on Control Systems Tech.}, \emph{arXiv preprint arXiv:2008.02379}, 2021.
	
\end{thebibliography}
\end{document}